\def\Bl{\rm{ Bl}}
\def\IC{{\bf C}}
\def\TheMagstep{\magstep1}	
\def\PaperSize{letter}		

\magnification=\magstep1

\let\:=\colon  
   \let\?=\overline

\let\Sum=\sum \def\sum{\Sum\nolimits}
\def\Bl{\rm{ Bl}}
\def\IC{{\bf C}} 
\def\IP{{\bf P}}

\def\and{\hbox{ and }}
	\def\Af{\hbox{\rm A$_f$}}

\def\DONE{*!*}
\def\NextDef #1 {\def\NextOne{#1}%
 \ifx\NextOne\DONE\let\next\relax
 \else\expandafter\xdef\csname#1\endcsname{\TheOp}
  \let\next\NextDef
 \fi \next}
\def\TheOp{\mathop{\rm\NextOne}}
 \NextDef 
  Projan Supp Proj Sym Spec Hom cod Ker dist
 *!*
\def\TheOp{{\cal\NextOne}}
\NextDef 
  E F G H I J M N O R S
 *!*
\def\TheOp{\hbox{\rm\NextOne}}
\NextDef 
 A ICIS 
 *!*
 
\def\item#1 {\par\indent\indent\indent\indent \hangindent4\parindent
 \llap{\rm (#1)\enspace}\ignorespaces}
 \def\inpart#1 {{\rm (#1)\enspace}\ignorespaces}
 \def\part {\par\inpart}

\catcode`\@=11		

\def\vfootnote#1{\insert\footins\bgroup
 \eightpoint 
 \interlinepenalty\interfootnotelinepenalty
  \splittopskip\ht\strutbox 
  \splitmaxdepth\dp\strutbox \floatingpenalty\@MM
  \leftskip\z@skip \rightskip\z@skip \spaceskip\z@skip \xspaceskip\z@skip
  \textindent{#1}\footstrut\futurelet\next\fo@t}

\def\p.{p.\penalty\@M \thinspace}
\def\pp.{pp.\penalty\@M \thinspace}
\def\(#1){{\rm(#1)}}\let\leftp=(
\def\activeleftp{\catcode`\(=\active}
{\activeleftp\gdef({\ifmmode\let\next=\leftp \else\let\next=\(\fi\next}}

\def\sct#1\par
  {\removelastskip\vskip0pt plus2\normalbaselineskip \penalty-250 
  \vskip0pt plus-2\normalbaselineskip \bigskip
  \centerline{\smc #1}\medskip}

\newcount\sctno \sctno=0
\def\sctn{\advance\sctno by 1 
 \sct\number\sctno.\quad\ignorespaces}

\def\dno#1${\eqno\hbox{\rm(\number\sctno.#1)}$}
\def\Cs#1){\unskip~{\rm(\number\sctno.#1)}}

\def\proclaim#1 #2 {\medbreak
  {\bf#1 (\number\sctno.#2)}\enspace\bgroup\activeleftp
\it}
\def\endproclaim{\par\egroup\medskip}
\def\pf{\endproclaim{\bf Proof.}\enspace}
\def\lem{\proclaim Lemma } 
\def\cor{\proclaim Corollary }	\def\thm{\proclaim Theorem }
\def\rmk#1 {\medbreak {\bf Remark (\number\sctno.#1)}\enspace}
\def\eg#1 {\medbreak {\bf Example (\number\sctno.#1)}\enspace}

\parskip=0pt plus 1.75pt \parindent10pt
\hsize29pc
\vsize44pc
\abovedisplayskip6pt plus6pt minus2pt
\belowdisplayskip6pt plus6pt minus3pt

\def\TRUE{TRUE}	
\ifx\DoublepageOutput\TRUE \def\TheMagstep{\magstep0} \fi
\mag=\TheMagstep

\newskip\vadjustskip \vadjustskip=0.5\normalbaselineskip
\def\centertext
 {\hoffset=\pgwidth \advance\hoffset-\hsize
  \advance\hoffset-2truein \divide\hoffset by 2\relax
  \voffset=\pgheight \advance\voffset-\vsize
  \advance\voffset-2truein \divide\voffset by 2\relax
  \advance\voffset\vadjustskip
 }
\newdimen\pgwidth\newdimen\pgheight
\def\letter{letter}\def\AFour{AFour}
\ifx\PaperSize\letter
 \pgwidth=8.5truein \pgheight=11truein 
 \message{- Got a paper size of letter.  }\centertext 
\fi
\ifx\PaperSize\AFour
 \pgwidth=210truemm \pgheight=297truemm 
 \message{- Got a paper size of AFour.  }\centertext
\fi

 \newdimen\fullhsize \newbox\leftcolumn
 \def\fulline{\hbox to \fullhsize}
\def\doublepageoutput
{\let\lr=L
 \output={\if L\lr
          \global\setbox\leftcolumn=\columnbox \global\let\lr=R%
        \else \doubleformat \global\let\lr=L\fi
        \ifnum\outputpenalty>-20000 \else\dosupereject\fi}%
 \def\doubleformat{\shipout\vbox{%
        \fulline{\hfil\hfil\box\leftcolumn\hfil\columnbox\hfil\hfil}%
				}%
		  }%
 \def\columnbox{\vbox
   {\makeheadline\pagebody\makefootline\advancepageno}%
   }%
 \fullhsize=\pgheight \hoffset=-1truein
 \voffset=\pgwidth \advance\voffset-\vsize
  \advance\voffset-2truein \divide\voffset by 2
  \advance\voffset\vadjustskip
 \let\firstheadline=\hfil
 
}
\ifx\DoublepageOutput\TRUE \doublepageoutput \fi

 \font\twelvebf=cmbx12		
 \font\smc=cmcsc10		

\def\eightpoint{\eightpointfonts
 \setbox\strutbox\hbox{\vrule height7\p@ depth2\p@ width\z@}%
 \eightpointparameters\eightpointfamilies
 \normalbaselines\rm
 }
\def\eightpointparameters{%
 \normalbaselineskip9\p@
 \abovedisplayskip9\p@ plus2.4\p@ minus6.2\p@
 \belowdisplayskip9\p@ plus2.4\p@ minus6.2\p@
 \abovedisplayshortskip\z@ plus2.4\p@
 \belowdisplayshortskip5.6\p@ plus2.4\p@ minus3.2\p@
 }
\newfam\smcfam
\def\eightpointfonts{%
 \font\eightrm=cmr8 \font\sixrm=cmr6
 \font\eightbf=cmbx8 \font\sixbf=cmbx6
 \font\eightit=cmti8 
 \font\eightsmc=cmcsc8
 \font\eighti=cmmi8 \font\sixi=cmmi6
 \font\eightsy=cmsy8 \font\sixsy=cmsy6
 \font\eightsl=cmsl8 \font\eighttt=cmtt8}
\def\eightpointfamilies{%
 \textfont\z@\eightrm \scriptfont\z@\sixrm  \scriptscriptfont\z@\fiverm
 \textfont\@ne\eighti \scriptfont\@ne\sixi  \scriptscriptfont\@ne\fivei
 \textfont\tw@\eightsy \scriptfont\tw@\sixsy \scriptscriptfont\tw@\fivesy
 \textfont\thr@@\tenex \scriptfont\thr@@\tenex\scriptscriptfont\thr@@\tenex
 \textfont\itfam\eightit	\def\it{\fam\itfam\eightit}%
 \textfont\slfam\eightsl	\def\sl{\fam\slfam\eightsl}%
 \textfont\ttfam\eighttt	\def\tt{\fam\ttfam\eighttt}%
 \textfont\smcfam\eightsmc	\def\smc{\fam\smcfam\eightsmc}%
 \textfont\bffam\eightbf \scriptfont\bffam\sixbf
   \scriptscriptfont\bffam\fivebf	\def\bf{\fam\bffam\eightbf}%
 \def\rm{\fam0\eightrm}%
 }

\def\today{\ifcase\month\or	
 January\or February\or March\or April\or May\or June\or
 July\or August\or September\or October\or November\or December\fi
 \space\number\day, \number\year}
\nopagenumbers
\headline={%
  \ifnum\pageno=1\firstheadline
  \else
    \ifodd\pageno\oddheadline
    \else\evenheadline\fi
  \fi
}
\let\firstheadline\hfill
\def\oddheadline{\eightpoint \rlap{\today}
 \hfil\headtitle\hfil\llap{\folio}}
\def\evenheadline{\eightpoint\rlap{\folio}
 \hfil\author\hfil\llap{\today}}
\def\headtitle{\title}

 \newcount\refno \refno=0	 \def\NoKey{*!*}
 \def\MakeKey{\advance\refno by 1 \expandafter\xdef
  \csname\TheKey\endcsname{{\number\refno}}\NextKey}
 \def\NextKey#1 {\def\TheKey{#1}\ifx\TheKey\NoKey\let\next\relax
  \else\let\next\MakeKey \fi \next}
 \def\RefKeys #1\endRefKeys{\expandafter\NextKey #1 *!* }
\def\SetRef#1 #2,#3\par{%
 \hang\llap{[\csname#1\endcsname]\enspace}%
  \ignorespaces{\smc #2,}
  \ignorespaces#3\unskip.\endgraf
 }
 \newbox\keybox \setbox\keybox=\hbox{[8]\enspace}
 \newdimen\keyindent \keyindent=\wd\keybox
\def\references{
  \bgroup   \frenchspacing   \eightpoint
   \parindent=\keyindent  \parskip=\smallskipamount
   \everypar={\SetRef}}
\def\endreferences{\egroup}

 \def\serial#1#2{\expandafter\def\csname#1\endcsname ##1 ##2 ##3
  {\unskip\ #2 {\bf##1} (##2), ##3}}
 \serial{ajm}{Amer. J. Math.}
  \serial {aif} {Ann. Inst. Fourier}
 \serial{asens}{Ann. Scient. \'Ec. Norm. Sup.}
 \serial{comp}{Compositio Math.}
 \serial{conm}{Contemp. Math.}
 \serial{crasp}{C. R. Acad. Sci. Paris}
 \serial{dlnpam}{Dekker Lecture Notes in Pure and Applied Math.}
 \serial{faa}{Funct. Anal. Appl.}
 \serial{invent}{Invent. Math.}
 \serial{ma}{Math. Ann.}
 \serial{mpcps}{Math. Proc. Camb. Phil. Soc.}
 \serial{ja}{J. Algebra}
 \serial{splm}{Springer Lecture Notes in Math.}
 \serial{tams}{Trans. Amer. Math. Soc.}

\def\UThin{\penalty\@M \thinspace\ignorespaces}
\def\relaxnext@{\let\next\relax}
\def\cite#1{\relaxnext@
 \def\nextiii@##1,##2\end@{\unskip\space{\rm[\SetKey{##1},\let~=\UThin##2]}}%
 \in@,{#1}\ifin@\def\next{\nextiii@#1\end@}\else
 \def\next{{\rm[\SetKey{#1}]}}\fi\next}
\newif\ifin@
\def\in@#1#2{\def\in@@##1#1##2##3\in@@
 {\ifx\in@##2\in@false\else\in@true\fi}%
 \in@@#2#1\in@\in@@}
\def\SetKey#1{{\bf\csname#1\endcsname}}

\catcode`\@=12 

\def\title{ The L\^e numbers of the square of a function and their
applications}
\def\author{Javier Fern\'andez de Bobadilla and Terence Gaffney}
\RefKeys BLS BS1 BS2 BG  DG D   Bo1 Bo2 Bo3 G G1 Ki
LP L-R L-T Mac Mas1 Mas2 Ma Ma1 Ma2 P P1 Si Te Ti Za
\endRefKeys

\def\topstuff{\leavevmode
\bigskip\bigskip
\centerline{\twelvebf \title}
\bigskip
\centerline{\author}
\medskip\centerline{\today}
\bigskip\bigskip}
\topstuff
\sct Introduction

L\^e numbers were introduced by Massey in \cite{Mas1}, \cite{Mas2} with
the purpose of numerically
controlling the topological properties of families of non-isolated
hypersurface singularities and describing the topology associated
with a function with non-isolated singularities. They are a
generalization of the Milnor number for isolated hypersurface
singularities.

The L\^e numbers are well-defined even for functions, $f$, which are not
square-free; in this case the singular set
of the function, denoted $S(f)$, has a component which is a component of
$V(f)$. In this note we will describe the relation
between the L\^e numbers of $f$ and $f^2$, where $f$ is square-free, and we
will use this connection in two applications.

If $f$ is not square-free, then the L\^e numbers of $f$ have a
contribution from the polar multiplicities of the components of $S(f)$
which are also components of
$V(f)$. This contribution is described precisely in lemma 1.1, and is the
basis for our two applications.

This paper can be viewed as a step in understanding the L\^e numbers of the composite of a mapping with a function. Here we study the composite of an arbitrary square-free $f$ and $z^2$. In \cite {G1}, the second author studied the composite of a mapping $G$ defining an ICIS (isolated complete intersection singularity)  with a Morse function general with respect to the discriminant of $G$.

The first application is concerned with the extent to which the L\^e
numbers are invariant in a family of functions which satisfy some
equisingularity condition.

L\^e and Ramanujam proved in \cite{L-R}
that, except in the surface case a family of isolated hypersurface
singularities with constant Milnor
number has constant Milnor fibration (in a differentiable sense)
and constant ambient topological type. Massey showed in
\cite{Mas1},\cite{Mas2} that a family of
hypersurface singularities with critical set of codimension at least three in
the ambient space and constant L\^e
numbers has constant Milnor fibration in the differentiable sense. Massey
then asked (in \cite{Ma}) whether
the constancy of the L\^e numbers is strong enough to ensure the constancy
of the embedded topological
type.

The first author gave a negative answer to this question in \cite{Bo3} by
showing, via explicit
examples, that the constancy of the L\^e numbers does not even imply the
constancy of the topological
type of the abstract link (the intersection of the central fiber with the
boundary of the Milnor sphere
thought as an abstract space, i.e., without its embedding into the Milnor
sphere). However he showed
that, if the dimension of the hypersurfaces is at least three, the
constancy of the L\^e numbers ensures
the constancy of the homotopy type of the abstract link.

Conversely, the Milnor number is a topological invariant. Massey asked in~\cite{Mas1} whether the L\^e
numbers could be topological invariants.
One should note that in an strict sense the answer to this question is negative:
if $f$ is a square-free holomorphic germ then $f$ and $f^2$
have the same ambient
topological type but different L\^e numbers. Therefore it makes sense to
reformulate Massey's question
using stronger notions instead of ambient topological invariance.
Reasonable ones seem to be topological
right-equivalence invariance or ambient topological invariance along a
family.
In this paper we will show that L\^e numbers are not, in general, invariants
for the seemingly strongest
purely topological notion: {\it the L\^e numbers are not necessarily
constant in a family of
topologically right-equivalent functions}.

The way of showing this is as follows: we show
that if a function $f$ defines an isolated
singularity then the L\^e numbers of its square $f^2$, with respect to a
generic coordinate system,
determine, and are determined by Teissier's $\mu^*$-sequence.
Therefore if $f_t$ is topologically trivial but not a $\mu^*$ constant
family, then $f_t^2$ is the desired
counterexample (see the details in section 2). This looks a bit
surprising: the L\^e numbers are
not topological invariants using the strongest topological equivalence,
and they are not strong enough to
control the topology of the abstract link in a family, but on the other
hand
they can be used to characterize ambient Whitney equisingularity for
families of isolated
singularities. It would be interesting to determine more precisely the
topological content of the L\^e
numbers.

Finally we notice that using squares of reduced functions we could not find any family of 
topologically right-equivalent functions such
that the top L\^e-number (L\^e number corresponding to the highest dimensional L\^e cycle) 
is not constant. In fact we will note that if in any family of topologically right-equivalent squares 
of functions defining isolated singularities the top L\^e 
number is constant then Zariski's multiplicity conjecture is true for families of isolated 
singularities and conversely. This could lead one to think that the top L\^e number is invariant
in families of topologically right-equivalent functions. We show that this is false by a 
counterexample (which is not a family of squares of functions).

The idea of looking at the L\^e numbers of the square of a function has
another, completely different,
interesting application. Using them we can get a quick proof of a new
formula for the Euler obstruction for hypersurfaces.
We explain this connection below.

The Euler obstruction is an idea introduced by MacPherson (\cite{Mac}) as a
key step in developing the notion of the Chern class for singular spaces.
Recently,
Brasselet, L\^e, Seade \cite{BLS} found a formula for the Euler
obstruction using a general linear form $L$. Their formula is:
$${\rm Eu}_X(0)=\sum_i \chi(V_i\cap B_\varepsilon \cap L^{-1}(t_0)) \cdot
{\rm Eu}_X(V_i),$$
where $\ B_\varepsilon$ is a small ball around $0$ in $\IC^N$, $t_0
\in
\IC
\setminus \{ 0 \}$ is sufficiently near $\{ 0 \}$ and ${\rm Eu}_X(V_i)$ is the
Euler
obstruction of $X$ at any point of the stratum $V_i$, $\{V_i\}$ a Whitney
stratification.

We use the L\^e numbers to derive a formula for the Euler obstruction at
the origin
of a hypersurface singularity. Our formula seems dual to that of BLS, in
that
our formula is a relative
formula ; we work not on $X$ but with the function $f$ defining $X$; hence
the terms in our formula come from components of the
exceptional divisor of the blowup by $J(f)$, not with a Whitney
stratification of $X$. Further, the main operators, Euler number
and complex link switch roles in these two formulas.

Although our formula applies directly only to hypersurface singularities,
L\^e and Teissier showed that it has wider utility.
In their seminal paper
\cite {L-T} theorem 6.4.1, they showed that the Euler obstruction at the
origin for any equidimensional reduced analytic set $X^d$,
is the
same as that of the hypersurface obtained by projecting $X^d$ to
$\IC^{d+1}$ using a generic projection.

In the first section we relate the L\^e numbers of $f^2$ with the L\^e
numbers and the
multiplicities of the relative and absolute polar varieties of $f$. In the
case in
which $f$ defines an isolated singularity, our formula boils down to the
relation between
the L\^e numbers of $f^2$ and the $\mu^*$-sequence of $f$. After that we
derive the relations
between L\^e numbers and topology predicted above. In the second section
we
compute our formula for the Euler obstruction. In the third, we use the
formula
to compute some examples.

In particular we are able to compute the
Euler obstruction of any $X$ defined by a function which only has
transverse $D(q,p)$
singularities off the origin. Examples of this class are singularities of
finite
codimension with respect to certain types of ideals
(see \cite{Za}, Example 11.3 in \cite{Bo2}). These singularities have been
studied and
are interesting in many
respects (see \cite{Bo2}, \cite{Bo3}, \cite{Si} and references therein).
This example also shows how our formula and that of \cite{BLS} can be used
together in the study of hypersurface singularities which
satisfy some transversality
conditions except at the origin.

When this work was completed David Massey pointed out to us that our Theorem 2.1 also follows from 
 the general index formula of Brylinski, Dubson, and Kashiwara, coupled with some results from his papers. We also thank Massey for a careful reading of an earlier preprint version of our paper.

\sctn The L\^e numbers of the square of a function

\noindent Setup: Let $f:\IC^n,0\to \IC,0$, $X=f^{-1}(0)$. One can define
the L\^e numbers of $f$ at $0$ as follows. Blow up $\IC^n,0$ by the
jacobian
ideal of $f$ denoted
$J(f)$. The blowup $B_{J(f)}(\IC^n)$ is a subset of
$\IC^n\times\IP^{n-1}$, with exceptional divisor $[E]$. By abuse of
notation we call subsets of the
form
$\IC^n\times H$ hyperplanes if $H$ is a hyperplane on $\IP^{n-1}$ and
denote their cycle class by $[H]$. Similarly, we denote the  cycle class of the intersection of  $r$ general hyperplanes by $[H]^{r}$. The L\^e cycle of $f$ of dimension $k$,
$0\le
k\le n-1$, denoted $\Lambda^{k}(f,0)$, is
$\pi_*([E]\cdot [H]^{n-k-1})$, where  $[H]^{n-k-1}$ is general for  $[E]$, $\pi_*$ is the map on cycles induced
from the projection of the blowup to $\IC^n$.
If $m_i$ is the multiplicity of the $i$-th component of the underlying
set of $\Lambda^{k}(f,0)$, and $e_i$ is the degree of the component,
then 
the $k$-th L\^e number of $f$ at $0$, denoted $\lambda^k(f,0)$ is
$\sum_im_ie_i$. (Cf.\cite{Ma} p103.)

Note that the components $E_i$ of the exceptional divisor are the conormal
varieties of the projection of the components to $\IC^n$.
(Recall that the conormal variety of an analytic set $X$ consists of the
tangent hyperplanes to the smooth part of $X$, and their limits at
singular points.) This follows because $f$ satisfies the \Af
condition at the generic point $z$ of $\pi(E_i)$ for the pair
$(\IC^n,\pi(E_i))$, and a dimension count shows that the set of limits is
exactly the set of hyperplanes which contain
the tangent space to $\pi(E_i)$ at $z$.

This last remark leads to a relation between L\^e numbers and polar
varieties. In preparation for a description of this connection,
recall that the polar varieties of an analytic set $X\subset\IC^n$, are
defined by intersecting the conormal variety of $X$, denoted $C(X)$ with
hyperplanes on
$\IC^n\times\IP^{n-1}$ and projecting to $X$.

Now, the L\^e cycles of
$f$ can be decomposed into two types of cycles-- the fixed cycles, whose
underlying sets are the images of the components of the
exceptional divisor of
$\Bl_{J(f)} (\IC^n)$, and the moving cycles
whose underlying sets are the the polar varieties of the images of the
components of the exceptional divisor. (The fixed cycles are the part of
the
L\^e cycle which is independent of the choice of hyperplane; the moving
cycles are the part that moves with choice of hyperplane.) The underlying
sets of the moving cycles
are the the polar varieties of the images of the components of the
exceptional divisor because the components of the exceptional divisor are
the conormals of the images. (Cf. Theorem 7.5 of \cite {Ma2} for a discussion of this point in the more general context of characteristic cycles of perverse sheaves.)

Label the $j$-th component of the underlying sets of the fixed cycles of
dimension $i$ by $|\Lambda^i_{j,F}(f)|$, and
 a generic point of the component by $z_{i,j}$. 
Denote the absolute polar variety of dimension $i$
of a complex analytic variety $X$ at the
origin by $\Gamma^i(X)$, and denote
the relative polar variety of $f$ of dimension
$i+1$ at the point $z$ by $\Gamma^{i+1}(f,z)$.

\lem 1 Suppose $f:\IC^n,0\to \IC,0$, $f$ is square-free, $X=f^{-1}(0)$.
Then
$$\lambda^i(f^2)=2\lambda^i(f)+m(\Gamma^i(X))+\Sum_{ k\ge
i,j}m(\Gamma^{k+1}(f,z_{k,j}))m(\Gamma^i(|\Lambda^k_{j,F}(f)|),0),$$
\noindent where $z_{k,j}$ is a generic point of $|\Lambda^k_{j,F}(f)|$.

\pf The underlying set of the exceptional divisor of $B_{J(f^2)}(\IC^n)$
is the same as that of $B_{J(f)}(\IC^n)$ with the addition of
the strict transform of $X$ as a component. The strict transform of $X$ is
just the conormal variety of $X$. This follows because the limits
of the tangent hyperplanes to the fibers of $f^2$ at a point of $X$ are
generically the tangent plane to $X$
at the point, since \Af holds generically.
Thus, the underlying sets of the fixed cycles of $f^2$ consist of $X$ and
those of $f$, while
the underlying sets of the moving cycles of $f^2$ consist of the polar
varieties of $X$ and the underlying sets of the moving cycles of $f$.

Since $f$ is square free, at a generic point of any component of $X$,
$J(f^2)=I(X)$,  hence the contribution to the L\^e cycle of dimension
$i$ from $X$ is just $m(\Gamma^i(X))$, the multiplicity of $\Gamma^i(X)$ at
$0$.

If $V^k$ is the underlying set of a fixed L\^e cycle of $f$, then to
calculate the degree of $V$ for $f^2$, work at the generic point $z$ of
$V$, pick a generic plane $P$ of dimension complementary to $V$, find the
relative polar curve of $f^2$ restricted to $P$, and calculate the degree
of $J(f^2)=fJ(f)$ on the curve. This degree is just the degree of $J(f)$
on the curve plus the degree of $f$ on the curve. The polar curve in this
case is just the
polar curve of $f$ restricted to $P$, so the first degree is the same as
the degree of $V$ for $J(f)$, while the second is the degree of $V$ plus
the
multiplicity of the polar curve (\cite{Ma} p. 24 proposition 1.23). Since
the multiplicity of the polar curve is just the multiplicity of
$\Gamma^{k+1}(f,z)$,
summing over all the fixed and moving components of $f$ of dimension $i$
we get $2\lambda^i(f)+
\Sum_{ k\ge i,j}m(\Gamma^{k+1}(f,z_j))m(\Gamma^i|\Lambda^k_{j,F}(f)|,0)$,
from which the result follows.

Notice that if $i=0$, then the formula becomes
$$\lambda^0(f^2)=2\lambda^0(f)+m(\Gamma^1(f,0)).$$
This is because there are no polar varieties of dimension $0$, so the only
component of $E_{J(f^2)}$
which contributes is that which projects to the origin.

Suppose now that $f$ defines an isolated singularity. Then the L\^e cycle
$\Lambda^i(f)$ is empty for $i>0$. Hence there is a unique fixed L\^e
cycle, which is of degree $0$ and is equal to $\mu(f)[0]$, where $\mu(f)$
is the Milnor number of $f$ and $[0]$ is the $0$-cycle given by the origin
of
$\IC^n$. In this context the formula given in Lemma 1.1 boils down to
$$\lambda^i(f^2)=m(\Gamma^i(X)),$$
if $i>0$ and to
$$\lambda^0(f^2)=\mu(f)+m(\Gamma^1(f,0)).$$

Notice that, using the fact that $f$ defines an isolated singularity, the
above equalities can be easily obtained as well from calculating the L\^e
cycles
with the procedure described in page 18 of \cite{Ma}.

Assume that the coordinate system $(x_1,...,x_n)$ with which we are
working is generic. Denote by $H_i$ the plane section defined by
$H_i:=V(x_1,...,x_i)$ for $i>0$, and $H_0:=\IC^n$. Define $f_i:=f|_{H_i}$.
Teissier's $\mu^*$-sequence associated to $f$
is the sequence of numbers defined by $\mu^{n}(f)=\mu(f)$ and
$\mu^{n-i}(f)=\mu(f_i)$, for $1\leq i\leq n-1$. Observe that
$\mu^1(f)=m(X)-1$.

The polar varieties of $X$ in these coordinates are given by
$$\Gamma^{n-1}(X)=X,$$
$$\Gamma^i(X)=V(f,\partial f/\partial x_{i+2},...,\partial f/\partial
x_n),$$
for $1\leq i\leq n-2$.

By the genericity of the coordinate system the set
$$V(x_1,...,x_i,f,\partial f/\partial x_{i+2},...,\partial f/\partial
x_n)$$
\noindent is isolated at the origin, and moreover we have the equalities
$$m(\Gamma^i(X))=I_0(H_i,\Gamma^i(X))=l(O_{\IC^n,0}/(x_1,...,x_k,f,\partial
f/\partial x_{i+2},...,\partial f/\partial x_n)).$$

On the other hand, as $f_i$ defines an isolated singularity, we have
$$\Gamma^{1}(f_i,0)=V(x_1,...,x_i,\partial f/\partial x_{i+2},...,\partial
f/\partial x_n),$$
for $0\leq i\leq n-1$. Hence
$$l(O_{\IC^n,0}/(x_1,...,x_k,f,\partial f/\partial
x_{i+2},...,\partial f/\partial x_n))=I_0(V(f),\Gamma^{1}(f_i,0)).$$
By \cite{LP}, for $i\leq n-2$, we have
$$I_0(V(f),\Gamma^{1}(f_i,0))=\mu(f_i)+\mu(f_{i+1})=\mu^{n-i}(f)+\mu^{n-i-1}(f).$$
We have also
$$m(\Gamma^1(f,0))=\mu^n(f)+\mu^{n-1}(f).$$

Combining the above equalities we have proved:

\lem 2 If $f$ defines an isolated singularity
$$\lambda_{n-1}(f^2)=m(X),$$
$$\lambda_i(f^2) =\mu^{n-i}+\mu^{n-i-1},$$
for $1\leq i\leq n-2$, and
$$\lambda_0(f^2)=2\mu^n(f)+\mu^{n-1}(f).$$
\endproclaim

By \cite{Te} and \cite{BS1} and Lemma 1.2 we immediately deduce:

\cor 3 Let $f_t$ be a family of isolated singularities holomorphically
depending on a parameter $t$. Then $f_t$ is ambient Whitney equisingular
if and only if the L\^e numbers of $f_t^2$ are independent of $t$.
\endproclaim

The family $f_t:=x_3^5+tx_2^6x_3+x_2^7x_1+x_1^{15}+x_4^2$ is a stabilisation
of the Brian\c con and Speder
example of $\mu$-constant but not $\mu^*$ constant family (see
\cite{BS2}). By \cite{Ti},
there exist positive $\epsilon$ and $\delta$ and a continuous family of
homeomorphisms
$$h_t:B_\epsilon\cap f_0^{-1}(D_\delta)\to B_\epsilon\cap
f_t^{-1}(D_\delta),$$
such that $f_t$ composed with $h_t$ gives $f_0$.

As $f_t^{-1}(D_\delta)=(f_t^2)^{-1}(D_{\delta^2})$ and $f_t^2$ composed with $h_t$ 
gives $f_0^2$ we have

\eg 4 The family $f^2_t$ is topologically right-equivalent trivial but its
L\^e numbers are not constant.

We note that we have not found examples of topologically right-equivalent functions which are squares 
of functions defining isolated singularities such that the
L\^e number corresponding to the L\^e cycle of maximal (the top L\^e number) 
dimension is different for them. The difficulty is explained by the following observation: 

\lem 5 Zariski's multiplicity conjecture for isolated singularities is equivalent to the following statement: let $f$ and $g$ be
function germs defining isolated singularities, if $f^2$ and $g^2$ are topologically right
equivalent then their top L\^e numbers are equal.
     
\pf Let us deduce first Zariski's multiplicity conjecture from the proposed statement. 
Suppose that $f$ and $g$  
have the same ambient topological type. According to \cite{Ki} they are topologically 
right-equivalent. By the argument above then $f^2$ and $g^2$ are also topologically right equivalent, 
and hence they have the same top L\^e number. As the top L\^e number of the square of a reduced 
function $f$ equals the multiplicity at the origin of $f$ we are done.

Suppose now that Zariski's multiplicity conjecture is true. If $f^2$ and $g^2$ are topologically 
right-equivalent then clearly $f$ and $g$ have the same embedded topological type. By Zariski's 
multiplicity conjecture they have then the same multiplicity, and hence $f^2$ and $g^2$ have the same
top L\^e number. We are done.

We describe now a counterexample to the topological right-equivalence invariance of the top L\^e 
number in general (for non-necessarily squares). Define germs $f,g_t:(\IC^3,0)\to\IC$ by 
$$f(x,y,z):=x^{15}+y^{10}+z^6$$
$$g_t(x,y,z):=xy+tz.$$
Note that if the $x,y,z,t$ have weights $(2,3,5,0)$ respectively then $f$ and $g_t$ are weighted homogeneous with weights $(30,5)$.
Further, the family of pairs $(f,g_t)$ defines a family of i.c.i.s. with constant Milnor number and 
non-constant multiplicity (the multiplicity of the i.c.i.s. is 12 for $t=0$ and $10$ for $t\neq 0$). 
The family is due to Henry and appears in \cite{BG}. Define $F_t:(\IC^3,0)\to\IC$ by
$$F_t:=f^2-g_t^{12}=(f+g^6_t)(f-g^6_t).$$
Here the form of $F_t$ is chosen so that $F_t$ is weighted homogeneous with weight $60$. It is easy to check that the singular set of the central fibre of $F_t$ is the intersection of the 
hypersurfaces $V(f+g^6_t)$ and $V(f-g^6_t)$, which is precisely $V(f,g_t)$. 
The generic transversal singularity is $A_{11}$ for any value of $t$. The first L\^e number of $F_t$ 
is the product of the multiplicity of $V(f,g_t)$ by the Milnor number of the generic transversal 
singularity. Hence the top L\^e number of the family $F_t$ is not constant. However:

\lem 6 The family $F_t$ is topologically right-equivalent trivial.

\pf We will prove this using the techniques of \cite{DG}. The basis of these techniques as described on  p341 of \cite{DG} applies 
to any function whether 
the singularities are non-isolated or not. We will describe how to construct a time dependent vector field $\xi(t,z)$ which is continuous 
everywhere and integrable, 
such that $DF_t(-1,\xi)=0$. Since the flow of $<-1,\xi>$ is C$^0$, this will show that the family is right-equivalent trivial.

To construct $\xi(t,z)$ we want to understand $J_z(F_t)$, the ideal generated by the partials of $F_t$ with respect to $(x,y,z)$.
It is helpful to think of $F$ as the composition of $X^2-Y^{12}$ and $(f,g_t)$. We can then view $J_z(F_t)$ as the image of the composite 
of two maps,  $D(X^2-Y^{12})\circ (f,g_t) $
and $[M_t]=[D(f_t,g_t),k_t]$ where $k_t$ generates the kernel of $D(X^2-Y^{12})\circ (f,g_t)$. It then follows that $h_12f_t-h_212g^{11}_t$ is in
$J_z(F_t)$ if and only if $(h_1,h_2)$ is in the $\O_4$ module $M_t$ generated by the columns of $[M_t]$. If $t=0$, then an easy bare hands computation shows that 
the $\O_3$ module $M_0$ has rank 2 except at the origin, hence contains $m^k_3\O^2_3$ for some $k$. Since in $[M_t]$, we are deforming the last two columns
by terms of weight greater than or equal to the weights of the entries in the last two columns, 
by Nakayama's lemma we can show that $M_t$ contains $m^k_3\O^2_4$.

Now we construct the control function and take the next steps in the construction of $\xi(t,z)$. Start with $x^{a/2}, y^{a/3}, z^{a/5}$ for $a$ very
 large and divisible by $30$, multiply each 
by $(F_t)_t/(12b)(xy-tz)^{11}$, which is the preimage of $ (F_t)_t$ , write this as  $M(v_x), M(v_y), M(v_z)$, where
 we can take $v_x,v_y,v_z$ weighted homogeneous vector fields. Now multiply the first equation 
$$x^{a/2}  (F_t)_t/(12b)(xy-tz)^{11}=M(v_x)$$

\noindent by the complex conjugate of $x^{a/2}$, do a similar thing to the other equations, sum the equations; the control 
function $\rho$ is the sum of products of the $x^{a/2}, y^{a/3}, z^{a/5}$ with their conjugates.

Let  $v'_x$ be the field made up of the first three components of $( \bar x)^{a/2}v_x$. Then the desired vector field is $1/\rho(v'_x+ v'_y+ v'_z)$ off the 
$t$ axis and $0$ on the $t$ axis.

It is clear that $\xi(t,z)$ is defined everywhere and real analytic off the $t$ axis. The continuity of $\xi(t,z)$ and the integrability of
$<-1,\xi(t,z)>$ follow from the material in \cite{DG} p341-343.

The technique used in this example covers a broad class of weighted homogeneous, 
non-isolated singularities, but we do not state a general result here.

It would be nice to answer:
 
{\bf Question 1.7}. For which functions is the top L\^e number an invariant of topological right 
equivalence?

\sctn The formula for the Euler obstruction

We continue with the notations of the previous section. At a generic point
$z_{i,j}$ of
$|\Lambda^i_{j,F}(f)|$, chose a generic plane $H$ of complementary
dimension to $i$. The complex link of $X\cap H$ at $z_{i,j}$
has the homotopy type of a bouquet of spheres, let
$b(|\Lambda^i_{j,F}(f)|,X)$ denote the number of spheres in this bouquet.
Denote the Euler obstruction to an analytic set $V$ at $x$ by ${\rm
Eu}(V,x)$. Given a complex analytic set, let $\chi(L,X,0)$ denote the
Euler
characteristic of the complex link of $X$ at $0$.

\thm 1 Suppose $f:\IC^n,0\to \IC,0$, $f$ is square-free, $X=f^{-1}(0)$.
Then
$${\rm
Eu}(X,0)=\chi(L,X,0)+\Sum_{i>0,j}(-1)^{n-i}b(|\Lambda^i_{j,F}(f)|,X){\rm
Eu}(|\Lambda^i_{j,F}(f)|,0).$$

\pf We know that the alternating sum of the L\^e numbers of $f^2$ is the
reduced Euler characteristic $\tilde \chi(M(f^2))$ of the Milnor fiber of
$f^2$
(\cite{Ma1}, Theorems II.10.3
and IV.3.5). The Milnor fiber
of $f^2$ is just two disjoint copies of the Milnor fiber of $f$, so the
reduced Euler characteristic of $f^2$ is:
$$\tilde \chi(M(f^2))=2\tilde \chi(M(f)+1.$$

Taking the alternating sum of the L\^e numbers of $f^2$ we get:

$$2\tilde \chi(M(f)+1=\sum_{i=0}^{n-1}(-1)^{n-1-i} \lambda^i(f^2)=$$
$$\sum_{i=0}^{n-1}(-1)^{n-1-i} (2\lambda^i(f)+m(\Gamma^i(X))+\Sum_{ k\ge
i,j}m(\Gamma^{k+1}(f,z_{k,j}))m(\Gamma^i(|\Lambda^k_{j,F}(f)|),0)) $$

Now $\sum_{i=0}^{n-1}(-1)^{n-1-i} (2\lambda^i(f))=2\tilde \chi(M(f))$, and
the alternating sum of the polar multiplicities
of the germ of an equidimensional analytic space $(X,x)$ is just ${\rm
Eu}(X,x)$ by the main theorem of \cite{L-T} (p476). Using these two facts
we get:

$$1={\rm Eu}(X,0)+(-1)^{n-1}m((\Gamma^1(f,0)))$$
$$+\Sum_{i>0,j}(-1)^{n-i-1}m(\Gamma^{k+1}(f,z_j)){\rm
Eu}(|\Lambda^i_{j,F}(f)|,0).$$

Solving for ${\rm Eu}(X,0)$ we get:

$${\rm Eu}(X,0)=1+(-1)^nm((\Gamma^1(f,0))$$
$$+\Sum_{i>0,j}(-1)^{n-i}m(\Gamma^{k+1}(f,z_{k,j})){\rm
Eu}(|\Lambda^i_{j,F}(f)|,0).$$

By a theorem of Massey's (p.365, \cite{Ma2}) we have that the multiplicity
of the relative polar variety of $f$ of dimension $k+1$ at the point $z_j$
is the $b(|\Lambda^k_{j,F}(f)|,X)$. Then
$\chi(L,X,0)=1+(-1)^{n-2}m(\Gamma^1(f))=1+(-1)^{n}m(\Gamma^1(f))$. So,
$${\rm Eu}(X,0)=\chi(L,X,0)+\Sum_{i>0,j}(-1)^{n-i}b(|\Lambda^i_{j,F}(f)|,X){\rm
Eu}(|\Lambda^i_{j,F}(f)|,0).$$

\sctn Examples

In the first two examples we recover the calculation of Eu$(X)$ in two
basic examples.

\eg 1 Suppose $f$ has an isolated singularity, then there are no
$\Lambda^i_{j,F}$ terms, and we get
${\rm Eu}(X,0)=\chi(L,X,0)=1+(-1)^{n-2}\mu^{n-1}(X,0)$.\vskip .1in

\eg 2 Suppose the singular locus of $f$, $S(f)$ consists of curves. Then
each curve is a fixed component of the L\^e cycles; choosing a
complementary
and generic hyperplane and intersecting with $X$ we are again in the
situation of an isolated hypersurface singularity, so
$b(|\Lambda^1_{j,F}(f)|,X)=\mu^{n-2}(X\cap H,z_j)=(e_j-1)(-1)^{n-3}$,
where $e_j$ is the value of the Euler obstruction to $X$ at a generic
point of $|\Lambda^1_{j,F}(f)|$. Now, ${\rm
Eu}(|\Lambda^1_{j,F}(f)|,0)=m(|\Lambda^1_{j,F}(f)|,0)$, by the theorem of
L\^e and Teissier,
so we recover the result of Dubson (\cite{D}),
$${\rm
Eu}(X,0)=\chi(L,X,0)+\Sum_j(-1)^{n-1}(-1)^{n-3}(e_j-1)=\chi(L,X,0)-\Sum_j(1-e_j).$$\vskip
.1in

In the next two examples we look at those $X$ which are defined by
functions which off the origin have only $A(k)$ singularities (example 3)
and $D(q,p)$ singularities (example 4). These ideas were originally
defined by Siersma and his students and are interesting because these
germs play
the same role for non-isolated hypersurface singularities
that finitely determined functions play for isolated hypersurface
singularities.

\vskip .1in

\eg 3 We say $f$ defined on $\IC^n$ has a singularity of type $A(d)$ at
$z$ if the singular set at $z$
has dimension $d$ and the germ of $f$ at $z$ has normal form
$z^2_1+\dots+z^2_{n-d}$.
Suppose $f$ is a function such that off the origin
$f$ has only
$A(d)$ singularities; denote the singular set by
$V$, then $V$ has an isolated singularity at the origin. Then
the fixed components of dimension $>0$ are exactly the components $V_i$ of $V$,
and $b(V_i,X)$ is $1$ for all components.
Meanwhile, since $V$ has an isolated singularity, ${\rm
Eu}(V_i,0)=\chi(L_{V_i},0),$ where
$\chi(L_{V_i},0)$ denotes the Euler characteristic of the complex link of
$V_i$ at $0$.
So ${\rm Eu}(X,0)=\chi(L,X,0)+(-1)^{n-d}\Sum_i\chi(L_{V_i},0).$

\vskip .1in

In preparation for the next example we discuss the $D(q,p)$ singularities.
These are singularities of functions in which the singular locus is
smooth.
If $f$ has a non-isolated singularity of type $D(p(p+1)/2,p)$,
then it is known that by a change of coordinates, $f$ has the normal form

$$f({\bf x},{\bf y})=(\Sum_{i\le j}^p
x_{i,j}y_iy_j)+y^2_{p+1}+\dots+y^2_{p+k}=[{\bf y}]^t[{\bf X}][{\bf
y}]+y^2_{p+1}+\dots+y^2_{p+k}.$$

Here $[{\bf X}]$ is a symmetric matrix with diagonal entries $x_{i,i}$
and off diagonal entries $1/2x_{i,j}$ and $n$, the dimension of the domain
of $f$
is $p+k+p(p+1)/2$.
In the notation $D(q,p)$, $p$ refers to the size of the matrix $[{\bf
X}]$ while the number of generators of the ideal $I=({\bf y})$ which
defines the singular locus of $f$
is $p+k$, while $q$ is the dimension of the singular set.
The smallest $q$ can be is $p(p+1)/2$. If $q>p(p+1)/2$, then the
additional coordinates do not appear in the normal form. (Cf. \cite{P},
\cite{P1}.)

\vskip .1in
\eg 4 Now, suppose $S(f)$ is an analytic set $V^d$ such that off the
origin $f$ has only $A(d)$ singularities or $D(q,p)$ singularities
appearing transversely
(see \cite{Bo2}, section 11.3).
Then, there are two types of fixed cycles of dimension $>0$ ; they are the
underlying sets of the components of $V^d$ and $D$,
the set of points where $f$ has a $D(q,p)$ singularity. (Assume $D$ has
dimension at least 1, otherwise
we are back in the last case.) The assumption on the behavior of $f$ off
the origin again implies that $V^d$ has an isolated singularity.
At a generic point $z_i$
of
$D_i$, a component of $D$, $f$ has the singularity type of a suspended
Whitney umbrella. So, $b(D_i,X)$ is the multiplicity of the relative
polar curve of $f$ at $z$, which is 1. So

$${\rm
Eu}(X,0)=\chi(L,X,0)+(-1)^{n-d}\Sum_i\chi(L_{V_i},0)+(-1)^{n-d-1}\Sum_i{\rm
Eu}(D_i,0).$$

Using our formula alone, at this point we would compute ${\rm Eu}(D_i,0)$
using the formula of L\^e and Teissier.
Notice that this is no more complicated than calculating the contributions
to the L\^e numbers coming from the $D_i$.

Combining our formula with that of Brasselet, L\^e and Seade we can make
further progress. By the transversality condition which $f$ satisfies,
a Whitney stratification of $D_i$ is given by the stratification by
$D(q,p)$ type. The Euler invariant at a generic point of each stratum
can be calclulated because the multiplicities of the polar varieties are
known. For the set of symmetric $p\times p$ square matrices
the polar multiplicities at the zero matrix of the set of matrices of
kernel rank 1 are given by:

$$m^{p(p+1)/2-i-1}=2^i{{p}\choose{p-i-1}},0\le i<p$$

$$m^{p(p+1)/2-i-1}=0,i\ge p.$$
(See \cite {G} for a proof.)

Then, at the generic point $z$ of each stratum of $D$, slicing by a
complementary generic plane $P$, the germ $X\cap P,z$ is
isomorphic to the germ of the set of symmetric $p\times p$ square matrices
of kernel rank 1 at the zero matrix where $p$ is determined from the
local nornmal form of $f$. Thus, the polar multiplicities
of $X\cap P,z$ at $z$ are the same as for the square matrices, and the
local Euler obstruction is the same as well.
Denote the value of the Euler obstruction to the germ of the set of
symmetric $p\times p$ square matrices
of kernel rank 1 at the zero matrix by ${\rm Eu}(p)$.

Looking at the formulas for the polar multiplicities we notice the
following relations:

$$2 {\rm Eu}(p)=2{{p}\choose{p-1}}+\dots+(-1)^{p-1}2^p{{p}\choose 0}$$

$$(-1)^{p-1} 2 {\rm Eu}(p)=2^p+\dots+(-1)^{p-1}2{{p}\choose{p-1}}$$

$$(-1)^{p-1} 2 {\rm
Eu}(p)+(-1)^p=2^p+\dots+(-1)^{p-1}2{{p}\choose{p-1}}+(-1)^p=(2-1)^p=1.$$

Hence $${\rm Eu}(p)={{1-(-1)^p}\over{(-1)^{p-1} 2}}.$$

Thus ${\rm Eu}(p)=0$ for $p$ even and ${\rm Eu}(p)=1$ for $p$ odd.

Denote the stratum of type $p$ of $D_i$ by $D_{i,p}$, by $D_p$ the union
of $D_{i,p}$ over $i$. Applying the formula of Brasselet, L\^e and Seade
to ${\rm Eu}(D,0)$, we get:

$${\rm Eu}(D_i,0)=\sum_{ p {\rm odd}} \chi(D_{i,p}\cap B_\varepsilon \cap
L^{-1}(t_0)). $$

Then we get:

$${\rm Eu}(X,0)=\chi(L,X,0)+(-1)^{n-d}\chi(L_{V},0)$$
$$+(-1)^{n-d-1}\sum_{ p {\rm odd}} \chi(D_{p}\cap B_\varepsilon \cap
L^{-1}(t_0)).$$

Based on the parity of $n-d$ this can be further simplified. Denote the
smooth points of $X$ by $X_0$. If $n-d$ is odd,

$${\rm Eu}(X,0)=\chi(X_0\cap B_\varepsilon \cap L^{-1}(t_0),0)+\sum_{ p
{\rm odd}} \chi(D_{p}\cap B_\varepsilon \cap L^{-1}(t_0)).$$

If $n-d$ is even, the story is a little more complicated. Let $A$ denote
the points on $S(X)$ where the germ of $f$ is of type $A(n-d)$.
Comparing the result of the BLS formula and ours we get:

$${\rm Eu}(X,0)=\chi(X_0\cap B_\varepsilon \cap L^{-1}(t_0),0)+
2\chi(A\cap B_\varepsilon \cap L^{-1}(t_0),0)$$
$$+\sum_{ p {\rm odd}} \chi(D_{p}\cap B_\varepsilon \cap L^{-1}(t_0))$$
and
$$\sum_{ p {\rm even}} \chi(D_{p}\cap B_\varepsilon \cap
L^{-1}(t_0))=0.$$

The procedure used to calculate this example should work whenever $f$
satisfies a set of transversality conditions off the origin.
From the transversality conditions it should be possible to see what the
fixed L\^e cycles are, and what the complex link of $X$ at a generic point
is.
Then calculate the Euler obstruction for each fixed cycle at the origin,
using the BLS formula. Here we only need a Whitney stratification of the
fixed cycle, not of $X$, and these should come from the transversality
conditions and the ``general situation".
The calculation of the Euler obstruction
should also
follow from the general situation as it does in this example.

\sct References

\references

BLS
J.P. Brasselet, L\^e D. T., J. Seade, { Euler obstruction and indices of
vector fields}, Topology
vol. 39 , 1193-1208, 2000

BS1
J.  Brian\c con, J.P. Speder, {\it Familles \' equisinguli\` eres de surfaces \` a singularit\'e isol\'ee},
C. R. Acad. Sci. Paris S\'er. A-B  280  (1975), Aii, A1013--A1016

BS2
J.  Brian\c con, J.P. Speder, {\it La trivialite topologique n'implique pas les conditions de Whitney},
C. R. Acad. Sci. Paris S\'er. A-B  280  (1975), no. 6, Aiii, A365--A367

BG
 R.-O. Buchweitz, G.-M. Greuel, {\it The Milnor number and deformations of complex curve singularities},  Invent. Math.  58  (1980), no. 3, 241--281

DG
J. Damon, T. Gaffney, {\it  Topological triviality of deformations of functions and Newton filtrations}, Invent. Math.  72  (1983), no. 3, 335--358

D
 	A. Dubson, {\it Classes caract\'eristiques des vari\'et\'es singuli\`eres}, C.R. Acad\'emie des Sciences Paris S\'erie A 287 (1978) 237--240

Bo1
J. Fern\'andez de Bobadilla, {\it Approximations of non-isolated 
singularities of finite codimension with respect to an isolated complete intersection 
singularity}.  Bull. London Math. Soc.  35  (2003),  no. 6, 812--816

Bo2
J. Fern\'andez de Bobadilla, {\it Relative morsification theory.}  
Topology  43  (2004),  no. 4, 925--982

Bo3
J. Fern\'andez de Bobadilla, {\it Answers to equisingularity questions.}
To appear in Invent. Math

G
T.Gaffney, {\it Invariants of $D(q,p)$ singularities}, Math CV/0508047

G1
T. Gaffney, { \it The multiplicity of pairs of modules and hypersurface 
singularities}, to appear  in Real and Complex Singularities (Luminy, 2004)

Ki
H. C. King, {\it Topological type of isolated critical points},
  Ann. Math. (2)  107  (1978), no. 2, 385--397 

LP
D. T. L\^e, B. Perron, {\it Sur la fibre de Milnor d'une singularit\'e isol\'ee en dimension
 complexe trois}, C. R. Acad. Sci. Paris S\'er. A-B  289  (1979), no. 2, A115--A118

L-R
D. T. L\^e and C. P. Ramanujam,
{\it The invariance of Milnor's number implies the invariance of the
topological type,}
\ajm 98 1976 67--78

L-T
 D. T. L\^e and B. Teissier, {\it Vari\'et\'es polaires locales et classes de Chern des vari\'et\'es singuli\`eres,}   Ann. of Math. (2)  114  (1981), no. 3,
457--491

Mac
R. D. MacPherson,
{\it Chern classes for singular algebraic varieties,} Ann. of Math. (2)
100 (1974), 423--432

Mas1
D. B. Massey, {\it The L\^e varieties. I.}  Invent. Math.  99  (1991),  no. 1, 357--376

Mas2
Massey, David B. {\it The L\^e varieties. II.}  Invent. Math.  104  (1991),  no. 1, 113--148

Ma
D. Massey {\it L\^e Cycles and Hypersurface Singularities,}
Springer Lecture Notes in Mathematics 1615 , (1995)

Ma1
D. Massey {\it Numerical Control over Complex Analytic Singularities}
Memoirs of the AMS,\#778, AMS 2003

Ma2
D. Massey, {\it Numerical invariants of perverse sheaves}, Duke Math. J.
73 (1994), no. 2, 307--369

P
R. Pellikaan,{\it Hypersurface singularities and
resolutions
of Jacobi modules,} Thesis, Rijkuniversiteit Utrecht, 1985

P1
R. Pellikaan, {\it Finite determinacy of functions with non-isolated
singularities,}
Proc. London Math. Soc. vol. 57, pp. 1-26, 1988

Si
D. Siersma, {\it The vanishing topology of non isolated singularities}.  
New developments in singularity theory (Cambridge, 2000),
447--472, NATO Sci. Ser. II Math. Phys. Chem., 21, Kluwer Acad. Publ., Dordrecht, 2001

Te
B. Teissier,
 {\it Cycles \'evanescents sections planes et conditions de Whitney,}
 in ``Singularit\'es \'a Carg\`ese," Ast\'erisque {\bf 7--8}, 1973

Ti
J.G. Timourian, {\it The invariance of Milnor's number implies topological triviality.}
Amer. J. Math.  99  (1977), no. 2, 437--446

Za
A. Zaharia,{\it Topological properties of certain singularities with critical locus
 a $2$-dimen-sional complete intersection.}  Topology Appl.  60  (1994),  no. 2, 153--171

\endreferences

\bye